\documentclass[12pt]{article}
\usepackage{amssymb}
\usepackage{amsthm,amstext,amscd,latexsym,amsmath,amsfonts,amssymb}
\renewcommand{\Box}{\square}
\theoremstyle{plain}
\newtheorem{theorem}{Theorem}
\newtheorem{lemma}[theorem]{Lemma}

\theoremstyle{definition}
\newtheorem{definition}[theorem]{Definition}

\newtheorem{example}{\sc Example}
\newtheorem{comment}[theorem]{Comment}
\theoremstyle{remark}
\newtheorem{remark}{\sc Remark}
\renewcommand{\Box}{\square}

\begin{document}
\title{\textbf{Regular obstructed categories and TQFT}}
\author{Steven Duplij\thanks{%
On leave of absence from Kharkov National University, Kharkov 61001, Ukraine}
\thanks{%
E-mail: Steven.A.Duplij@univer.kharkov.ua} \thanks{%
Internet: http://gluon.physik.uni-kl.de/\~{}duplij}\hskip0.7em and Wladyslaw
Marcinek\thanks{%
E-mail: wmar@ift.univ.wroc.pl} \\
Institute of Theoretical Physics, University of Wroclaw,\\
Pl. Maxa Borna 9, 50-204 Wroc{l}aw, Poland}
\date{}
\maketitle

\begin{abstract}
A proposal of the concept of $n$-regular obstructed categories is given. The
corresponding regularity conditions for mappings, morphisms and related structures in
categories are considered. An n-regular TQFT is introduced. It is shown the connection
of time reversibility with the regularity.
\end{abstract}

PACS. 03. 65. Fd - Algebraic methods in quantum theory

\newpage

\section{Introduction}
In the generalized histories approach \cite{gell} to quantum theory the whole universe
is represented by a class of 'histories'. In this approach the standard Hamiltonian
time-evolution is replaced by a partial semigroup called a 'temporal support'. A
possible realization of such program can be described in terms of cobordism manifolds
and corresponding categories \cite{isham}. The temporal support arises naturally as a
cobordism $M$, where the boundary $\partial M$ of $M$ is a disjoint sum of the
'incoming' boundary manifold $\Sigma_0$ and the 'outgoing' one $\Sigma_1$. This means
that the cobordism $M$ represents certain quantum process transforming $\Sigma_0$
into $\Sigma_1$. In other words, $\Sigma_1$ is a time consequence of $\Sigma_0$.
Obviously, we have two opposite possibilities to declare which boundary is the
initial one.

Let $N$ be a cobordism with the 'outgoing' boundary of $M$ as its 'incoming boundary'
and $\Sigma_2$ as the 'outgoing boundary'. Then there is a cobordism $N\circ M$ whose
incoming boundary is $\Sigma_0$, and the outgoing one is $\Sigma_3$. In this case we
say that these two cobordisms are glued along $\Sigma_1$. Such gluing of cobordisms
up to diffeomorphisms define a partial semigroup operation. One can consider cobordism
with several incoming and outgoing boundary manifolds. The class of possible
histories can be represented by gluing of cobordisms in several different ways. Hence
there is the corresponding coherence problem for such description.

Let $Cob$ be a category of cobordisms, where the boundary $\partial M$ of $M\in Cob$
is a disjoint sum of the 'incoming' boundary manifold $\Sigma_0$ and the 'outgoing'
one $\Sigma_1$. There is also the cylinder cobordism $\Sigma\times[0,1]$ such that
$\partial(\Sigma\times[0,1])=\Sigma\amalg\Sigma^{\ast}$. The class of boundary
components is denoted by $Cob_0$. According to Atiyah \cite{at}, Baez and Dolan
\cite{bad}, the TQFT is a functor $\mathcal{F}$ from the category $Cob$ to the
category $Vect$ of a finite-dimensional vector spaces. This means that $\mathcal{F}$
sends every manifold $\Sigma\in Cob_0$ into vector space $\mathcal{F}(\Sigma)$ such
that
\begin{equation}
\begin{array}{ccc}
\mathcal{F}(\Sigma^{\ast})=(\mathcal{F}(\Sigma))^{\ast}, & \mathcal{F}%
(\Sigma_0\amalg\Sigma_1)=(\mathcal{F}\Sigma_0)\otimes(\mathcal{F}\Sigma_1), &
\mathcal{F}(\emptyset)=I,
\end{array}
\end{equation}
and a cobordism $M(\Sigma_0,\Sigma_1)$ to a mapping $\Phi(M)\in
lin_{I}(\mathcal{F}\Sigma_0, \mathcal{F}\Sigma_1)$ such that $\mathcal{F} (\Sigma
\times[0,1])=id_{\mathcal{F}\Sigma}$, where $I$ is a field, and $\Sigma^{\ast}$ is the
same manifold $\Sigma$ but with the opposite orientation. Kerler \cite{ker} found
examples of categories formed by some classes of cobordism manifolds preserving some
operations like the disjoint sum or surgery. It was discussed by Baez and Dolan
\cite{bad} that it is not easy to describe such categories in a coherent way. Crane
\cite{cra,cran} applied the category theory to an algebraic structure of the quantum
gravity.

The idea of regularity as generalized inverse was firstly introduced by by von Neumann
\cite{neu} and applied by Penrose for matrices \cite{pen1}. Let $\mathbb{R}$ be a
ring. If for an element $a\in\mathbb{R}$ there is an element $a^{\star}$ such that
\begin{equation}
\begin{array}{cc}
aa^{\star }a=a, & a^{\star }aa^{\star }=a^{\star }, \label{reg}
\end{array}
\end{equation}
then $a$ is said to be regular and $a^{\star}$ is called a generalized inverse of
$a$. Generalizing transition from inverses to regularity is a widely used method of
abstract extension of various algebraic structures. The intensive study of such
regularity and related directions was developed in many different fields, e.g.
generalized inverses theory \cite{rabson,rao/mit,nashed}, semigroup theory
\cite{mun/pen,cli5,lal1,howie,lawson}, and \cite{dup6,dup11}, supermanifold theory
\cite{duplij,dup18,dup-hab}, Yang-Baxter equation in endomorphism semigroup and
braided almost bialgebras \cite{fangli3,fangli2,fangli1}, weak bialgebras, week Hopf
algebras \cite{nill}, category theory \cite{dav/rob}.

In this paper we are going to study certain class of categories which can be useful
for the study of quantum histories with non-reversible time, quantum, gravity and
field theory. The regularity concept for linear mappings and morphisms in categories
are studied. Higher order regularity conditions are described. Commutative diagrams
are replaced by `semicommutative' ones. The distinction between commutative and
`semicommutative' cases is measured by a non-zero obstruction proportional to the
difference of some self-mappings $e^{\left( n\right) }$ from the identity. This
allows to 'regulaarize' the notion of categories, functors and related algebraic
structures. It is interesting that this procedure is unique up to an equivalence
defined by invertible morphisms. Our regularity concept is nontrivial for equivalence
classes of noninvertible morphisms. The regular version of TQFT is a natural
application of the presented here formalism. In this case the n-regularity means that
a time evolution is nonreversible, although repeated after n steps, but up to a
classes of obstructions. Our considerations are based on the concepts of generalized
inverse \cite{nashed,dav/rob}, and semisupermanifolds \cite{duplij}.

The paper is organized as follows. In the Section II we consider linear mappings
without requirement of `invertibility'. If $f:X\rightarrow Y$ is a linear mapping,
then Instead of the inverse mapping $f^{-1}:Y\rightarrow X$ we use less restricted
`regular' $f^{\star}$ one by extending `invertibility' to `regularity' according to
the following relations
\begin{equation}
f\circ f^{\star}\circ f  = f,  \quad f^{\star}\circ f\circ f^{\star} = f^{\star }.
\label{2r}
\end{equation}
We also propose some higher regularity conditions. In Section III the higher
regularity notion is extended to morphisms of categories. Commutative diagrams are
replaced by semicommutative ones. The concept of regular cocycles of morphisms in an
category is described. An existence theorem for these cocycles is given. The
corresponding generalization of certain categorical structures as tensor operation,
algebras and coalgebras etc... to our higher regularity case is given in the Section
IV. Regular equivalence classes of cobordism manifolds and the corresponding
structures are considered in the Section V. An $n$-regular TQFT is introduced as an
$n$-regular obstructed category represented some special classes of cobordisms called
'interactions'. Our study is not complete, it is only a proposal for new algebraic
structures related to topological quantum theories.

\section{Generalized invertibility and regularity}

Let $X$ and $Y$ be two linear spaces over a field $k$. We use the following
notation. Denote by $id_{X}$ and $id_{Y}$ the identity mappings $%
id_{X}:X\rightarrow X$ and $id_{Y}:Y\rightarrow Y$. If $f:X\rightarrow Y$is a linear
mapping, then the image of $f$ is denoted by $Im f$, and the kernel by $Ker f$.

We are going to study here some generalizations of the standard concept of
invertibility properties of mappings. Our considerations are based on the article of
Nashed \cite{nashed}. Let $f:X\rightarrow Y$ be a linear mapping. If $f\circ
f^{-1}_r=id_{Y}$ for some $f^{-1}_r:Y\rightarrow X$, then $f$ is called a
\textit{retraction}, and $f^{-1}_r$ is the right
inverse. Similarly, if $f^{-1}_l\circ f=id_{X}$, then it is called a \textit{%
coretraction}, $f^{-1}_l$ is the left inverse of $f$. A mapping $f^{-1}$ is
called an inverse of f if and only if it is both, right and left inverse of $%
f$.

This standard concept of inverses is in many cases to strong to be fulfilled. To
obtain more weak conditions one has to introduced the following `regularity'
conditions
\begin{equation}
f\circ f_{in}^{\star }\circ f=f,  \label{3r}
\end{equation}
where $f_{in}^{\star}:Y\longrightarrow X$ is called an \textit{inner inverse}%
, and such $f$ is called \textit{regular}. Similar ``reflexive regularity'' conditions
\begin{equation}
f_{out}^{\star }\circ f\circ f_{out}^{\star }=f_{out}^{\star }  \label{4r}
\end{equation}
defines an \textit{outer inverse} $f_{out}^{\star }$. Notice that in general $%
f_{in}^{\star }\neq f_{out}^{\star }\neq f^{-1}$ or it can be that $f^{-1}$ does not
exist at all.

\begin{definition}
A mapping $f$ satisfying one of the condition (\ref{3r}) or (\ref{4r}) is said to be
regular or 2-regular. A generalized inverse of a mapping $f$ is a mapping $f^{\star }$
which is both inner and outer inverse $f^{\star }=f_{in}^{\star }=f_{out}^{\star }$.
\end{definition}

\begin{lemma}
If $f_{in}^{\star }$ is an inner inverse of $f$, then a generalized inverse $%
f^{\star }$ exists, but need not be unique.
\end{lemma}

\textbf{Proof:} If $f_{in}^{\star }$ is an inner inverse, then
\begin{equation}
f^{\star }=f_{in}^{\star }\circ f\circ f_{in}^{\star }  \label{fin}
\end{equation}
is always both inner and outer inverse i. e. generalized inverse. It follows from
(\ref{fin}) that both regularity conditions (\ref{3r}) and (\ref{4r}) hold.

\hfill $\Box$\newline

\begin{definition}
Let us define two operators $\mathcal{P}_{f}:Y\rightarrow Y$ and $\mathcal{P}%
_{f^{\star }}:X\rightarrow X$ by relations
\begin{equation}
\mathcal{P}_{f}:=f\circ f^{\star },\quad \mathcal{P}_{f^{\star }}:=f^{\star }\circ f,
\end{equation}
\end{definition}

\begin{lemma}
These operators satisfy
\begin{equation}
\begin{array}{ll}
\mathcal{P}_{f}\circ \mathcal{P}_{f}=\mathcal{P}_{f}, & \mathcal{P}_{f}\circ
f=f\circ \mathcal{P}_{f^{\star }}=f \\
\mathcal{P}_{f^{\star }}\circ \mathcal{P}_{f^{\star }}=\mathcal{P}_{f^{\star }},
& \mathcal{P}_{f^{\star }}\circ f^{\star }=f^{\star }\circ \mathcal{P}%
_{f}=f^{\star }.
\end{array}
\end{equation}
\end{lemma}

\hfill $\Box$\newline

\begin{lemma}
If $f^{\star }$ is the generalized inverse of the mapping $f$, then the following
properties are obvious
\begin{equation}
\begin{array}{ll}
Imf=Im(f\circ f^{\star }), & Ker(f\circ f^{\star })=Kerf^{\star }, \\
Im(f^{\star }\circ f)=Imf^{\star }, & Ker(f^{\star }\circ f)=Kerf.
\end{array}
\end{equation}
In addition there are two decompositions
\begin{equation}
\begin{array}{ll}
X=Imf^{\star }\bigoplus Kerf, & Y=Imf\bigoplus Kerf^{\star },
\end{array}
\end{equation}
The restriction $f\mid _{Imf^{\star }}:Imf^{\star }\rightarrow Imf$ is one to one
mapping, and operators $P_{f},P_{f^{\star }}$ are projectors of $Y, X$ onto $Im f, Im
f^{\star }$, respectively.
\end{lemma}
\hfill $\Box$\newline

\begin{theorem}
Let $f:X\rightarrow Y$ be a linear mapping. If $P$ and $Q$ are projectors
corresponding to the following two decompositions
\begin{equation}
\begin{array}{ll}
X=M\bigoplus Kerf, & Y=Imf\bigoplus N,
\end{array}
\end{equation}
respectively, then exist unique generalized inverse of $f$, and
\begin{equation}
f^{\star }:=i\circ \tilde{f}^{-1}\circ Q,
\end{equation}
where $\tilde{f}:=f\mid _{M}$, and $i:M\hookrightarrow X$.
\end{theorem}
\hfill $\Box$\newline Here we try to construct higher analogs of generalized inverses
and regularity conditions (\ref{3r})--(\ref{4r}). Let us consider two mappings
$f:X\rightarrow Y$ and $f^{\star}:Y\rightarrow X$ and
introduce two additional mappings $f^{\star\star}:X\rightarrow Y$ and $%
f^{\star\star\star}:Y\rightarrow X$. We propose here the following higher regularity
condition
\begin{equation}
\begin{array}{l}
f\circ f^{\star}\circ f^{\star\star}\circ f^{\star\star\star}\circ f=f,
\end{array}
\end{equation}
This equation define a 4-regularity condition. By cyclic permutations we obtain
\begin{equation}
\begin{array}{l}
f^{\star}\circ f^{\star\star}\circ f^{\star\star\star}\circ f\circ
f^{\star}=f^{\star}, \\
f^{\star\star}\circ f^{\star\star\star}\circ f\circ f^{\star}\circ
f^{\star\star}=f^{\star\star}, \\
f^{\star\star\star}\circ f\circ f^{\star}\circ f^{\star\star}\circ f^{\star\star\star}
=f^{\star\star\star}. \label{3f}
\end{array}
\end{equation}
By recursive considerations we can propose the following formula of $n$%
-regularity
\begin{equation}
\begin{array}{l}
\underbrace{f\circ f^{\star}\circ f^{\star\star}\ldots\circ
f^{\overset{n=2k-1}{\overbrace {\star\star\ldots\star}}}\circ f}_{n+1}=f,
\label{nreg1}
\end{array}
\end{equation}
where $n=2k$, $k=1, 2, \ldots$ and their cyclic permutations. Observe that for a not
unique $\overset{2k-1}{\overbrace{\star\star\ldots\star}}$-operation the following
formula
\begin{align}
\left( g\circ f\right)^{\overset{2k-1}{\overbrace{\star\star\ldots\star}}} &
=f^{\overset{2k-1}{\overbrace{\star\star\ldots\star}}}\circ g^{\overset{2k-1}{%
\overbrace{\star\star\ldots\star}}}.  \label{gf1}
\end{align}
leads to a difficulty. If the above operation is defined up to an equivalence, then
the difficulty can be overcomes. We can introduce `higher projector' by the relation
\begin{equation}
\mathcal{P}_{f}^{\left( n\right) }=f\circ f^{\star}\circ f^{\star\star}\ldots\circ
f^{\overset{2k-1}{\overbrace{\star\star\ldots\star}}},\quad n=2k. \label{pf}
\end{equation}

\begin{lemma}
It is easy to check the following properties
\begin{align}
\mathcal{P}_{f}^{\left( n\right) }\circ f & =f.  \label{pff2}
\end{align}
and $\mathcal{P}_{f}^{\left(n\right) }\circ\mathcal{P}%
_{f}^{\left(n\right) }=\mathcal{P}_{f}^{\left(n\right) }, n=2k$.
\end{lemma}
\hfill $\Box$\newline
For a given $n=2k$ all $f^{\star}, f^{\star\star}, \ldots f^{%
\overset{n-1}{\overbrace{\star\star\ldots\star}}}$ are different, and, for instance,
$\left( f^{\star}\right) ^{\star}\neq f^{\star\star}$. The existence of analogous
conditions for odd $n$ is a problem.
\begin{theorem}
Let $f:X\rightarrow Y$ be a linear mapping. If $P$, and $Q$ are projectors
corresponding to the following two decompositions
\begin{equation}
\begin{array}{ll}
X=M\bigoplus Kerf, & Y=Imf\bigoplus N,
\end{array}
\end{equation}
respectively, and
\begin{equation}
\begin{array}{ll}
f^{\star }\mid _{Imf}=f^{\star \star \star }\mid _{Imf}, &
\end{array}
\end{equation}
then the 4-regularity condition of $f$ can reduced to the two 2-regularity conditions
\begin{equation}
\begin{array}{ll}
f\circ f^{\star }\circ f=f, & f^{\star }\circ f^{\star \star }\circ f^{\star
}=f^{\star }.
\end{array}
\end{equation}
\end{theorem}
\hfill $\Box$\newline

\section{Semicommutative diagrams and regular obstructed categories}

In the previous section we considered mappings and regularity properties for two
given spaces $X$ and $Y$, because we studied various types of inverses. Now we will
extend these consideration to any number of spaces and introduce semicommutative
diagrams (firstly introduced in \cite{duplij}).

A \textit{directed graph} $\mathfrak{C}$ is a pair $\{\mathfrak{C}_{0},
\mathfrak{C}_{1}\}$ and a pair of functions
\begin{equation}
\mathfrak{C}_{0}
\begin{array}{c}
s \\
\leftleftarrows  \\
t
\end{array}
\mathfrak{C}_{1}
\end{equation}
where elements of $\mathfrak{C}_{0}$ are said to be \textit{objects}, elements of
$\mathfrak{C}_{1}$ are said to be \textit{arrows} or \textit{morphisms}, $sf$ is said
to be a \textit{domain (or source)} of $f$, and $tf$ is a \textit{codomain (or
target)} of $f\in \mathfrak{C}_{1}$. If $sf= X\in\mathfrak{C}_o$, and $tf= Y\in
\mathfrak{C}_0$, then we use the following notation $X\stackrel{f}{\longrightarrow}Y$
and
\begin{equation}
\begin{array}{c}
\mathfrak{C}(X, Y):=\{f\in\mathfrak{C}_1 : sf=X, tf=Y\}.
\end{array}
\end{equation}
We denote by $End(X)$ the collection of all morphisms defined on $X$ into itself, i.
e. $End(X):=\mathfrak{C}(X, X)$, $X\in\mathfrak{C}_0$.

Two arrows $f, g\in\mathfrak{C}_1$ such that $tf= sg$ are said to be composable. If
in addition $sf=X$, $sg=tf= Y$, and $tg=Z$, then we use the notation
$X\stackrel{f}{\longrightarrow}Y\stackrel{g}{\longrightarrow}Z$. In this case a
composition $g\circ f$ of two arrows $f:X\rightarrow Y$ and $g:Y\rightarrow Z$ can be
defined as an arrow $X\stackrel{g\circ f}{\rightarrow}Z$. The associativity means that
$h\circ\left( g\circ f\right) =\left( h\circ g\right) \circ f=h\circ g\circ f$. An
identity $'id'$ in $\mathfrak{C}$ is an inclusion $X\in\mathfrak{C}_0\rightarrowtail
id_X\in End(X)$ such that
\begin{equation}
\begin{array}{c}
f\circ id_X = id_Y\circ f = f.
\end{array}
\end{equation}
for every $X, Y\in\mathfrak{C}_1$, and $X\stackrel{f}{\longrightarrow}Y$.

A directed graph $\mathfrak{C}$ equipped with associative composition of composable
arrows and identity satisfying some natural axioms is said to be a \textit{category}
\cite{mitchell,maclane1}. If $\mathfrak{C}$ is a category, then right cancellative
morphisms are \textit{epimorphisms} which satisfy $g_{1}\circ f=g_{2}\circ f
\Longrightarrow$ $g_{1}=g_{2}$, where $g_{1,.2}:Y\rightarrow Z$ and left cancellative
morphisms are \textit{monomorphisms} which satisfy $f\circ h_{1}=f\circ h_{2}
\Longrightarrow$ $h_{1}=h_{2}$, where $h_{1,.2}:Z\rightarrow X$. A morphisms
$X\stackrel{f}{\longrightarrow}Y$ is invertible means that there is a morphisms
$Y\stackrel{g}{\longrightarrow}X$ such that $f\circ g=id_Y$ and $g\circ f=id_X$.
Instead of such invertibility  we can use the regularity condition (\ref{3r}), i.e.
$f\circ g\circ f=f$, where $g$ plays the role of an inner inverse \cite{nashed}.
\medskip

\begin{center}
\setlength{\unitlength}{.3in} \begin{picture}(15,3.5) \put(4.6,3.2){{ $f$}}
\put(4.6,1.7){{ $g$}} \put(4,2.8){\vector(1,0){2.1}} \put(6.1,2.4){\vector(-1,0){2.1}}
\put(6.3,3.4){\small`Regularization'} \put(8,2.5){$\Longrightarrow$} \put(10.6,1.6){{
$g$}} \put(10.6,3.2){{ $f$}} \put(10,2.9){\vector(1,0){2.1}}
\put(12.1,2.3){\vector(-1,0){2.1}} \put(10,2.6){\vector(1,0){2.1}}
\put(0,2.5){\large$n=2$} \put(4,0.2){\makebox(1,1){\small Invertible morphisms}}
\put(12,0.2){\makebox(1,1){\small Noninvertible (regular) morphisms}}
\end{picture}
\end{center}

Usually, for three objects $X, Y, Z$ and three morphisms $f:X\rightarrow Y$ and
$g:Y\rightarrow Z$ and $h:Z\rightarrow X$ one can have the `invertible' triangle
commutative diagram $h\circ g\circ f=id_{X}$. Its regular extension has the form
\begin{equation}
f\circ h\circ g\circ f=f.  \label{fgh}
\end{equation}
Such a diagram

\begin{center}
\setlength{\unitlength}{.3in}
\begin{picture}(15,5)
\put(5.8,1.5){\vector(-1,1){1.68}} \put(4.6,4.2){{ $f$}} \put(3.5,2){{ $h$}}
\put(4,3.5){\vector(1,0){2.1}} \put(6.1,3.2){\vector(0,-1){2}}
\put(6.3,4.3){\small`Regularization'} \put(8,2.5){$\Longrightarrow$} \put(6.5,2){{
$g$}} \put(9.5,2){{ $h$}} \put(11.8,1.5){\vector(-1,1){1.68}} \put(10.6,4.2){{ $f$}}
\put(10,3.8){\vector(1,0){2.1}} \put(10,3.4){\vector(1,0){2.1}}
\put(12.1,3.2){\vector(0,-1){2}} \put(12.5,2){{ $g$}} \put(0,2.5){\large$n=3$}
\put(4,0){\makebox(1,1){\small Reversible morphisms}} \put(12,0){\makebox(1,1){\small
Noninvertible (regular) morphisms}}
\end{picture}
\end{center}
can be called a \textit{semicommutative diagram}. By cyclic permutations of
(\ref{fgh}) we obtain
\begin{equation}
\begin{array}{l}
h\circ g\circ f\circ h = h, \\
g\circ f\circ h\circ g = g. \label{2rf}
\end{array}
\end{equation}
These formulae define the concept of 3-regularity.

\begin{definition}
A mapping $f:X\longrightarrow Y$ satisfying conditions (\ref{fgh}) and (\ref {2rf})
is said to be 3-regular. The mapping $h:Z\longrightarrow X$ is called the first
3-inversion and the mapping $g:Y\longrightarrow Z$ the second one.
\end{definition}

The above concept can be expanded to any number of objects and morphisms.
\begin{definition}
Let $\mathfrak{C}=(\mathfrak{C}_{0}, \mathfrak{C}_{1})$ be a directed graph. An
\textit{n-regular cocycle} $(X, f)$ in $\mathfrak{C}$, $n=1, 2, \ldots$, is a sequence
of composable arrows in $\mathfrak{C}$
\begin{equation}
\begin{array}{c}
X_{1}\stackrel{f_{1}}{\longrightarrow }X_{2}\stackrel{f_{2}}{\longrightarrow }\cdots
\stackrel{f_{n-1}}{\longrightarrow }X_{n}\stackrel{f_{n}}{\longrightarrow }X_{1},
\label{regu}
\end{array}
\end{equation}
such that
\begin{equation}
\begin{array}{c}
f_{1}\circ f_{n}\circ \cdots \circ f_{2}\circ f_{1}=f_{1},\\
f_{2}\circ f_{1}\circ \cdots \circ f_{3}\circ f_{2}=f_{2},\\
\\
f_{n}\circ f_{n-1}\circ \cdots \circ f_{1}\circ f_{n}=f_{n}, \label{cyp}
\end{array}
\end{equation}
and
\begin{equation}
\begin{array}{c}
e_{X_{1}}^{\left( n\right) }:=f_{n}\circ \cdots \circ f_{2}\circ f_{1}\in
End(X_{1}),\\
e_{X_{2}}^{\left( n\right) }:=f_{1}\circ \cdots \circ f_{3}\circ f_{2}\in
End(X_{2}),\\
\\
e_{X_{n}}^{\left( n\right) }:=f_{n-1}\circ \cdots \circ f_{1}\circ f_{n}\in
End(X_{n}). \label{cocy}
\end{array}
\end{equation}
\end{definition}

\begin{definition}
Let $(X, f)$ be an n-regular cocycle in $\mathfrak{C}$, then the correspondence
$e_X^{(n)}:X_i\in\mathfrak{C}_0\mapsto e_{X_i}^{(n)}\in End(X_i)$, $i=1, 2, \ldots,
n$, is called an $n$-regular cocycle obstruction structure on $(X, f)$ in
$\mathfrak{C}$.
\end{definition}

\begin{lemma}
We have the following relations
\begin{equation}
\begin{array}{ccc}
f_i\circ e_{X_i}^{\left( n\right) }= f_{i},&e_{X_{i+1}}^{(i)}\circ f_i =
f_{i},&e_{X_i}^{\left( n\right) }\circ e_{X_i}^{\left( n\right) }=e_{X_i}^{\left(
n\right) }. \label{ob}
\end{array}
\end{equation}
for $i=1, 2, \ldots ,n (\mbox{mod}\;n)$.
\end{lemma}
{\bf Proof:} The lemma  simply follows from relations (\ref{cyp}) and (\ref{cocy}).
\hfill $\Box$\newline

\begin{definition}
An $n$-regular obstructed category is a directed graph $\mathfrak{C}$ with an
associative composition and such that every object is a component of an n-regular
cocycle.
\end{definition}

\begin{example}
If all obstruction are equal to the identity $e_{X_i}^{\left( n\right) }=id_{X_i}$,
and
\begin{equation}
\begin{array}{c}
f_{n}\circ \cdots \circ f_{2}\circ f_{1}=id_{X_1},\\
f_{1}\circ \cdots \circ f_{3}\circ f_{2}=id_{X_2},\\
\\
f_{n-1}\circ \cdots \circ f_{1}\circ f_{n}=id_{X_n}. \label{cocyt}
\end{array}
\end{equation}
then the sequence (\ref{regu}) is trivially n-regular. Observe that the trivial
$2$-regularity is just the usual invertibility, hence every grupoid $G$ is a trivially
$2$-regular obstructed category. We are interested with obstructed categories
equipped with some obstruction different from the identity.
\end{example}

\begin{definition}
The minimum number $n=n_{obstr}$ such that $e_{X}^{(n)}\neq id_X$ is called the
obstruction degree.
\end{definition}

\begin{example}
Every inverse semigroup $S$ is a nontrivial $2$-regular obstructed category. It has
only one object, morphisms are the elements of $S$.
\end{example}

\begin{theorem}
Let $\mathfrak{C}$ be a category, and
\begin{equation}
\begin{array}{c}
X_{1}\overset{f_{1}}{\longrightarrow }X_{2}\overset{f_{2}} {\longrightarrow }\cdots
\overset{f_{n-1}}{\longrightarrow}X_{n} \overset{f_{n}}{\longrightarrow
}X_{1}\label{seq}
\end{array}
\end{equation}
be a sequence of morphisms of category $\mathfrak{C}$. Assume that there is a sequence
\begin{equation}
\begin{array}{c}
Y_{1}\overset{\widetilde{f}_{1}}{\longrightarrow }Y_{2}\overset{\widetilde{f}_{2}}
{\longrightarrow }\cdots\overset{\widetilde{f}_{n-1}}{\longrightarrow}Y_{n}
\overset{\widetilde{f}_{n}}{\longrightarrow }Y_{1},\label{seqp}
\end{array}
\end{equation}
where $Y_i$ is a subobject of $X_i$ such that there is a collection of mappings
$\pi_i:X_i\rightarrow Y_i$ and $\iota:Y_i\rightarrow X_i$ satisfying the condition
$\pi_i\circ\iota_i = id_{Y_i}$ for $i=1, 2, \ldots, n$. If in addition
\begin{equation}
\begin{array}{l}
\widetilde{f}_n\circ\cdots\widetilde{f}_2\circ\widetilde{f}_1 = id_{Y_1},\\
\widetilde{f}_1\circ\cdots\widetilde{f}_3\circ\widetilde{f}_2 = id_{Y_2},\\
\ldots\\
\widetilde{f}_{n-1}\circ\cdots\widetilde{f}_{1}\circ\widetilde{f}_n = id_{Y_n},
\label{commd}
\end{array}
\end{equation}
and
\begin{equation}
\begin{array}{c}
f_{i} := \iota_{i+1}\circ\widetilde{f}_{i}\circ\pi_i \label{lift}
\end{array}
\end{equation}
then the sequence (\ref{seq}) is an $n$-regular cocycle.
\end{theorem}
{\bf Proof:} The corresponding obstruction structure is given by
\begin{equation}
\begin{array}{l}
e^{(n)}_{X_i} = \iota_i\circ\pi_i\label{obb}
\end{array}
\end{equation}
If $x\in Ker f_{1}$, then the theorem is trivial, if $x\in X_i\setminus Ker f_1$, then
we obtain
$$
\begin{array}{rl}
(f_1\circ f_n\circ\cdots\circ f_2\circ f_1)(x)& = \iota_2\circ\widetilde{f}_1
\circ\pi_1\circ\iota_1 \circ\widetilde{f}_n\circ\cdots\circ\widetilde{f}_2
\circ\widetilde{f}_1\circ\pi_1(x)\\& = \iota_{2}\circ\widetilde{f}_{1}\circ\pi_1 =
f_1(x),
\end{array}
$$
where the condition (\ref{commd}) and (\ref{lift}) has been used. We can calculate
all cyclic permutations in an similar way. \hfill $\Box$\newline

\begin{example}
There is an  $n$-regular obstructed category $\mathfrak{C}=(\mathfrak{C}_0,
\mathfrak{C}_1)$, where $\mathfrak{C}_0=\{X_i:i=1, \ldots, n\;(\mbox{mod} n) \}$ and
$\mathfrak{C}_1=\{f_i:i=1,\ldots, n\;(\mbox{mod} n)\}$ are described in the above
theorem.
\end{example}

\begin{definition}
Let $(X,f), (Y,g)$ be two $n$-regular cocycles in $\mathfrak{C}$. An $n$-regular
cocycle morphism $\alpha:(X,f)\rightarrow (Y,g)$ is a sequence of morphisms
$\alpha:=(\alpha_{1},\ldots,\alpha_{n})$ such that the diagram
\begin{equation}
\begin{array}{rccccccccl}
& X_{1} & \overset{f_{1}}{\longrightarrow} & X_{2} & \overset{f_{2}}{%
\longrightarrow}\cdots\overset{f_{n-1}}{\longrightarrow} & X_{n} & \overset{%
f_{n}}{\longrightarrow} & X_{1} &  &  \\
& \downarrow\alpha_{1} &  & \downarrow\alpha_{2} &  & \downarrow\alpha_{n} &
& \downarrow\alpha_{1} &  &  \\
& Y_{1} & \overset{g_{1}}{\longrightarrow} & Y_{2} & \overset{g_{2}}{%
\longrightarrow}\cdots\overset{g_{n-1}}{\longrightarrow} & Y_{n} & \overset{%
g_{n}}{\longrightarrow} & Y_{1} &  &
\end{array}
\end{equation}
is commutative. If every component $\alpha_i$ of $\alpha$ is invertible, then $\alpha$
is said to be an n-regular cocycle equivalence.
\end{definition}

It is obvious that the n-regular cocycle equivalence is an equivalence relation.

\begin{definition}
Let $\mathfrak{C}$ be an $n$-regular obstructed category. A collection of all
equivalence classes of n-regular cocycles in $\mathfrak{C}$ and corresponding
$n$-regular cocycle morphisms is denoted by $\Re eg^{(n)}(\mathfrak{C})$ and is said
to be an $n$-regularization of $\mathfrak{C}$.
\end{definition}

\begin{comment}
It is obvious that the n-regular cocycle equivalence is an equivalence relation.
Equivalence classes of this relation are just elements of $\Re
eg^{(n)}(\mathfrak{C})$. Our $n$-regular cocycles and obstruction structures are
unique up to an invertible n-regular cocycle morphisms. If $[(X, f)]$ is an
equivalence class of n-regular cocycles, then there is the corresponding class of
n-regular obstruction structures $e_X^{\left( n\right) }$ on it. The correspondence
is a one to one.
\end{comment}

\section{Regularization of functors and related structures}

We are going to introduce the concept of n-regular functors, natural transformations,
involution, duality, and so on. All our definitions of are in general case the same
like in the usual category theory \cite{maclane1}, but the preservation of the
identity $id_{X}$, is replaced by the requirement of preservation of obstructions
$e_X^{\left( n\right) }$ up to the n-regular cocycle equivalence.

It is known that for two usual categories $\mathfrak{C}$ and $\mathfrak{D}$ a functor
$\mathcal{F}:\mathfrak{C}\rightarrow\mathfrak{D}$ is defined as a pair of mappings
$(\mathcal{F}_0, \mathcal{F}_1)$, where $\mathcal{F}_0$ sends objects of
$\mathfrak{C}$ into objects of $\mathfrak{D}$, and $\mathcal{F}_1$ sends morphisms of
$\mathfrak{C}$ into morphisms of $\mathfrak{D}$
\begin{equation}
\begin{array}{c}
\mathcal{F}_1(f\circ g)=\mathcal{F}_1(f)\circ\mathcal{F}_1(g),\quad \mathcal{F}_1
id_X=id_{\mathcal{F}_0 X}.
\end{array}
\end{equation}
for $X\in \mathfrak{C}_{0}$, $\mathcal{F}X\in\mathcal{D}_{0}$.

Let $\mathfrak{C}$ and $\mathfrak{D}$ be two $n$-regular obstructed categories. We
postulate that all definitions are formulated on every n-regular cocycle $(X, f)$ in
$\mathfrak{C}$ up to the n-regular cocycle equivalence, and $i=1, 2,\ldots$(mod n).
\begin{definition}
An or $n$-regular cocycle functor
$\mathcal{F}^{(n)}:\mathfrak{C}\rightarrow\mathfrak{D}$ is a pair of mappings
$(\mathcal{F}^{(n)}_0, \mathcal{F}^{(n)}_1)$, where $\mathcal{F}^{(n)}_0$ sends
objects of $\mathfrak{C}$ into objects of $\mathfrak{D}$, and $\mathcal{F}^{(n)}_1$
sends morphisms of $\mathfrak{C}$ into morphisms of $\mathfrak{D}$ such that
\begin{equation}
\begin{array}{l}
\mathcal{F}_1^{(n)}(f_i\circ
f_{i+1})=\mathcal{F}_1^{(n)}(f_i)\circ\mathcal{F}^{(n)}_1(f_{i+1}), \quad
\mathcal{F}_1^{\left( n\right) }
\left( e_{X_i}^{\left( n\right) }\right) = e_{\mathcal{F}_0(X_i)}^{(n)}, \label%
{fe}
\end{array}
\end{equation}
where $X\in\mathfrak{C}_{0}$.
\end{definition}

\begin{lemma}
Let $\mathfrak{C}$ and $\mathfrak{D}$ be $n$-regular obstructed categories, and let
\begin{equation}
\begin{array}{c}
X_{1}\stackrel{f_{1}}{\longrightarrow }X_{2}\stackrel{f_{2}}{\longrightarrow
}\cdots \stackrel{f_{n-1}}{\longrightarrow }X_{n}\stackrel{f_{n}}{%
\longrightarrow }X_{1}\label{regul}
\end{array}
\end{equation}
be an $n$-regular cocycle in $\mathfrak{C}$. If $\mathcal{F}^{\left( n\right)
}:\mathfrak{C}\rightarrow \mathfrak{D}$ is $n$-regular cocycle functor, then
\begin{equation}
\begin{array}{c}
\mathcal{F}^{(n)}(f_i)\circ e^{(n)}_{X_i}=\mathcal{F}^{(n)}(f_i).
\end{array}
\end{equation}
\end{lemma}
\textbf{Proof:} It is a simple calculation
\begin{equation}
\mathcal{F}^{(n)}(f_i)=\mathcal{F}^{(n)}\left(f\circ e^{(n)}_{X_i}\right)=\mathcal{%
F}^{(n)}\left(f\right)\circ\mathcal{F}^{(n)}\left(e^{(n)}_{X_i}\right)=
\mathcal{F}^{(n)}(f_i)\circ e^{(n)}_{\mathcal{F}_0 X_i}.
\end{equation}
\hfill $\Box$\newline Multifuncors can be regularized in a similar way.

Let $\mathcal{F}^{(n)}$ and $\mathcal{G}^{(n)}$ be two $n$--regular cocycle morphisms
of the category $\mathfrak{C}$ into the category $\mathfrak{D}$.

\begin{definition}
An $n$-regular natural transformation $s:\mathcal{F}^{(n)}\rightarrow
\mathcal{G}^{(n)}$ of $\mathcal{F}^{(n)}$ into $\mathcal{G}^{(n)}$ is  a collection
of functors $s = \{s_{X_i}:\mathcal{F}_0(X_i)\rightarrow\mathcal{G}_0(X_i) \}$ such
that
\begin{equation}
s_{X_{i+1}}\circ \mathcal{F}^{(n)}_1(f_i) = \mathcal{G}^{(n)}_1(f_i)\circ s_{X_i},
\end{equation}
for $f_i:X_i\rightarrow X_{i+1}$.
\end{definition}

\begin{definition}
An $n$-regular obstructed monoidal category $\mathfrak{C}\equiv\mathfrak{C}(\otimes,
I)$ can be defined as usual, but we must remember that instead of the identity
$id_X\otimes id_Y = id_{X\otimes Y}$ we have an obstruction structure
$e_X^{(n)}=\{e_{X_i}^{\left( n\right) }\in End(X_i); n=1,2,...\}$ satisfying the
condition
\begin{equation}
e^{(n)}_{X_i\otimes Y_i} = e^{(n)}_{X_i}\otimes e^{(n)}_{Y_i}  \label{mul}
\end{equation}
\end{definition}
for every two n-regular cocycles $(X, f)$ and $(Y, f')$.

Let $\mathfrak{C}$ be an $n$-regular obstructed monoidal category. We introduce an
$\ast$-operation in $\mathfrak{C}$ as a function which send every object $X_i$ into
object $X^{\ast}_i$ called the dual of $X$,
\begin{equation}
\begin{array}{c}
X_i^{\ast\ast} = X_i, \quad(X_i\otimes Y_i)^{\ast} = X_i^{\ast}\otimes Y_i^{\ast},
\end{array}
\end{equation}
reverse all arrows
\begin{equation}
\begin{array}{c}
(f\circ g)^{\ast}=g^{\ast}\circ f^{\ast}.
\end{array}
\end{equation}
The category $\mathfrak{C}$ equipped with such $\ast$-operation is called an
\textit{$n$-regular obstructed monoidal category with duals}.
\begin{lemma}
Let $\mathfrak{C}$ be an n-regular obstructed monoidal category with duals. If $(X,
f)$ is an $n$-regular cocycle in $\mathfrak{C}$, then there is a corresponding
$n$-regular cocycle $(X^{\ast}, f^{\ast})$ in $\mathfrak{C}^{\ast}$, called the dual
of $(X, f)$.
\end{lemma}

\textbf{Proof:} If we reverse all arrows in $(X, f)$ and replace all objects by the
corresponding duals, then we obtain $(X^{\ast}, f^{\ast})$, where
\begin{equation}
\begin{array}{c}
X^{\ast}_{1}\stackrel{f^{\ast}_{n}}{\rightarrow}X^{\ast}_{n}\stackrel{%
f^{\ast}_{n-1}}{\rightarrow}\cdots\stackrel{f^{\ast}_{2}}{\rightarrow }%
X^{\ast}_{2}\stackrel{f^{\ast}_{1}}{\rightarrow }X^{\ast}_{1}\label{regub}
\end{array}
\end{equation}
is a sequence such that
\begin{equation}
\begin{array}{cc}
f^{\ast}_{1}\circ f^{\ast}_{n}\circ\cdots\circ f^{\ast}_{2}\circ
f^{\ast}_{1}=f^{\ast}_{1}, & e_{X^{\ast}_{1}}^{\left( n\right)
}:=f^{\ast}_{n}\circ\cdots\circ f^{\ast}_{2}\circ f^{\ast}_{1}.\label{egub}
\end{array}
\end{equation}
where $f^{\ast}_{i}:X^{\ast}_{i+1}\rightarrow X^{\ast}_i$, $i=1,\ldots, n$, and
$X^{\ast}_{n+1}\equiv X^{\ast}_1$ is the dual. We have corresponding relations for
all cyclic permutations. \hfill $\Box$

\begin{definition}
An n-regular pairing $g_{\mathfrak{C}}$ in an n-regular obstructed monoidal category
$\mathfrak{C}$ can be defined in an analogy to the usual case as a collection of
mappings
\begin{equation}
g_{\mathfrak{C}}=\{g_{X_i}\equiv\langle-|-\rangle_{X_i}:X_i^{\ast}\otimes
X_i\rightarrow I\}
\end{equation}
satisfying some natural consistency conditions and in addition the following
regularity relations
\begin{equation}
\begin{array}{c}
g_{X_{i+1}}\circ (f_i^{\ast}\otimes f_i) = g_{X_i},
\end{array}
\end{equation}
and
\begin{equation}
\begin{array}{l}
\langle e^{(n)}_{X^{\ast}_i}X^{\ast}_i\mid X_i\rangle_{X_i} = \langle X^{\ast}_i\mid
e^{(n)}_{X_i}X_i\rangle_{X_i}, \label{dual1ty}
\end{array}
\end{equation}
where $(X, f)$ is a regular $n$-cocycle in $\mathfrak{C}$, and let $(X^{\ast},
f^{\ast})$ is the corresponding duals.
\end{definition}

It is known that an associative algebra in an ordinary category is an object
$\mathcal{A}$ of this category such that there is an multiplication
$m:\mathcal{A}\otimes\mathcal{A}\rightarrow\mathcal{A}$ which is also a morphism of
this category satisfying some axioms like the associativity, the existence of the
unity.
\begin{definition}
Let $\mathfrak{C}$ be an $n$-regular obstructed monoidal category. An $n$-regular
cocycle algebra $\mathcal{A}$ in the category $\mathfrak{C}$ is an object of this
category equipped with an associative multiplication $m:\mathcal{A}\otimes\mathcal{A}
\rightarrow\mathcal{A}$ such that
\begin{equation}
m\circ(e^{(n)}_A\otimes e^{(n)}_A) = e^{(n)}_A\circ m . \label{regal}
\end{equation}
\end{definition}
Obviously such multiplication not need to be unique.

One can define an $n$-regular cocycle coalgebra or bialgebra in a similar way. A
comultiplication $\triangle :\mathcal{A}\longrightarrow\mathcal{A} \otimes
\mathcal{A}$ can be  regularized according to the relation
\begin{equation}
\triangle\circ e^{(n)}_{\mathcal{A}} = (e^{(n)}_{\mathcal{A}}\otimes
e^{(n)}_{\mathcal{A}})\circ\triangle . \label{regco}
\end{equation}
\begin{definition}
Let $\mathcal{A}$ be an $n$-regular cocycle algebra. If $\mathcal{A}$ is also regular
coalgebra such that $\Delta \left( ab\right) =\Delta \left( a\right) \Delta \left(
b\right) $, then it is said to be an $n$-regular cocycle almost bialgebra.
\end{definition}
If $\mathcal{A}$ is an $n$-regular cocycle algebra, then we denote by $hom_m
(\mathcal{A}, \mathcal{A} )$ the set of morphisms $s\in
hom_{\mathcal{C}}(\mathcal{A}, \mathcal{A} )$ satisfying the condition
\begin{equation}
s\circ m = m\circ (s\otimes s).
\end{equation}
Let $\mathcal{A}$ be an $n$-regular cocycle almost bialgebra. We define the
convolution product
\begin{equation}
s\star t := m\circ (s\otimes t)\circ\triangle,
\end{equation}
where $s, t\in hom_{m} (\mathcal{A}, \mathcal{A} )$. If $\mathcal{A}$ is a regular
$n$-cocycle almost bialgebra, then  the convolution product is regular.

\begin{definition}
An $2$-regular cocycle almost bialgebra $\mathcal{H}$ equipped with an element $S\in
hom_m (\mathcal{H},\mathcal{H})$ such that
\begin{equation}
\begin{array}{cc}
S\star id_{\mathcal{H}}\star S = S,&id_{\mathcal{H}}\star S\star id_{\mathcal{H}} =
id_{\mathcal{H}}.
\end{array}
\end{equation}
is said to be an $2$-regular cocycle almost Hopf algebra $\mathcal{H}$.
\end{definition}
The above definition is a regular analogy of week Hopf algebras considered in
\cite{nill}. Similar algebras has been also considered in \cite{fangdup}.

\begin{lemma}
If $\mathcal{A}$ is an $n$-regular cocycle algebra, then there is an $n$-regular
cocycle coalgebra $\mathcal{A}^{\ast}$ such that
\begin{equation}
\begin{array}{l}
\langle\triangle (\xi),x_1 \otimes x_2\rangle =\langle\xi,m(x_1\otimes x_2)\rangle,
\label{dual4}
\end{array}
\end{equation}
where $x_1, x_2\in\mathcal{A}, \xi\in\mathcal{A}^{\ast}$.
\end{lemma}
{\bf Proof:} Let us apply the regularity condition (\ref{regal}) to the above duality
condition (\ref{dual4}). Then the lemma follows from relations (\ref{mul}),
(\ref{regco}), and (\ref{dual1ty}). \hfill $\Box$\newline
\begin{lemma}
Let $\mathcal{A}$ be an $n$-regular cocycle almost bialgebra. Then the dual
$\mathcal{A}^{\ast}$ is also $n$-regular cocycle almost bialgebra
\begin{equation}
\begin{array}{l}
\langle\triangle (\xi),x_1 \otimes x_2\rangle
=\langle\xi,m(x_1\otimes x_2)\rangle,\\
\langle\widehat{m}(\xi\otimes\zeta) , x_1\otimes x_2\rangle =\langle\xi\otimes\zeta ,
\widehat{\triangle}x\rangle. \label{dual3}
\end{array}
\end{equation}
\end{lemma}
\hfill $\Box$\newline Let $\mathcal{A}$ be an $n$-regular cocycle algebra. Then we
can define a left $n$-regular cocycle $\mathcal{A}$-module as an object equipped with
a $\mathcal{A}$--module action $\rho_{M}:\mathcal{A}\otimes M\longrightarrow M$ such
that
\begin{equation}
\begin{array}{l}
\rho_{M}\circ (m\otimes id_M) = \rho_{M}\circ(id_{\mathcal{A}}\otimes\rho_{M}),\\
\rho_M\circ (e^{(n)}_{\mathcal{A}}\otimes e^{(n)}_M) = e^{(n)}_{M} \circ\rho_M .
\label{rmod}
\end{array}
\end{equation}
If $\mathcal{A}$ is an $n$-regular cocycle coalgebra, then one can define an n-regular
cocycle comodule $M$ in a similar way. For a coaction $\delta_M
:\mathcal{A}\rightarrow\mathcal{A}\otimes M$ of $\mathcal{A}$ on $M$ we have the
following regularity condition
\begin{equation}
\delta_M \circ (e^{(n)}_{\mathcal{A}}\otimes e^{(n)}_M) = e^{(n)}_{M}\circ\varrho_M ,
\label{rcom}
\end{equation}
\begin{remark}
Observe that we have the following duality between $\mathcal{A}$-module action
$\rho_{M}:\mathcal{A}\otimes M\longrightarrow M$ and $\mathcal{A}^{\ast}$-comodule
coactions $\delta_{M^{\ast}}:\mathcal{A}^{\ast}\rightarrow M^{\ast}\otimes
\mathcal{A}^{\ast}$
\begin{equation}
\begin{array}{l}
\langle\delta_{M^{\ast}}(\xi),a\otimes x\rangle =\langle\xi,\varrho_M(a\otimes
x)\rangle, \label{dual5}
\end{array}
\end{equation}
where $a\in\mathcal{\mathcal{A}}, x\in M, \xi\in\mathcal{A}^{\ast}$.
\end{remark}

\section{Regular cobordisms and TQFT}

Let $Cob$ be a directed graph of cobordisms whose objects $Cob_0$ are $d$-dimensional
compact smooth and oriented manifolds without boundary and whose arrows are classes
of cobordism manifolds with boundaries. We would like to discuss the corresponding
$n$-regular cocycles and their meaning. For this goal we use here a parametrization
such that the boundary $\partial M$ is a multiconnected space, a disjoint sum of the
'incoming' boundary manifold $\Sigma_{in}$ and the 'outgoing' one $\Sigma_{out}$. We
call them 'physical'. The empty boundary component is also admissible.  Let
$\Sigma_0, \Sigma_1\in Cob_0$, then the disjoint sum is denoted by
$\Sigma_0\amalg\Sigma_1$. For a manifold $\Sigma\in Cob_0$ there is the corresponding
manifold $\Sigma^{\ast}$ with the opposite orientation.

We wish to represent quantum processes of certain physical system by cobordism
manifolds $M$ with the 'incoming' boundary manifold $\Sigma_0$ (an 'input'),and the
'outgoing' one $\Sigma_0$, (an 'output'). The 'incoming' boundary manifold $\Sigma_0$
represents an initial conditions of the system, the 'outgoing' boundary represents
the final configuration, and the cobordism manifolds represent possible interaction
of the system. Note that the same cobordism manifold $M$ but with different boundary
parametrization represent different physical processes!
\begin{definition}
An 'interaction' is a triple $_{\Sigma_0}\mathcal{M}_{\Sigma_1}$, where the 'incoming'
boundary manifold $\Sigma_0$ is multiconnected space with $m$ components and the
'outgoing' one $\Sigma_1$ is equipped with $n$ components, and $\mathcal{M}$ is a
class of cobordism manifolds up to parametrization preserving diffeomorfisms,
$\Sigma_0, \Sigma_1\in Cob_0, \mathcal{M}\in Cob_1$.
\end{definition}

\begin{definition}
The 'opposite interaction' of $_{\Sigma_0}\mathcal{M}_{\Sigma_1}$ is the 'interaction'
$_{\Sigma_1}\mathcal{M}^{op}_{\Sigma_0}$ with reversed boundary parametrization, i. e.
the 'incoming' boundary of $\mathcal{M}$ is the 'outgoing' boundary of
$\mathcal{M}^{op}$ and vice versa,.
\end{definition}

\begin{example}
A 'collapsion' of $\Sigma\in Cob_0$ is an arbitrary 'interaction' of the forms
$_{\Sigma}\mathcal{M}_{\emptyset}$, this means the 'incoming' boundary is $\Sigma$
and the 'outgoing' boundary is empty. The corresponding 'expansion' of $\Sigma$ is
the opposite of the collapsion.
\end{example}
\begin{definition}
Let us denote by $\mathfrak{C}ob=(\mathfrak{C}ob_0, \mathfrak{C}ob_1)$ a directed
graph whose objects are $\mathfrak{C}ob_0\equiv Cob_0$ and arrows $\mathfrak{C}ob_1$
are 'interactions'. A composition  of two 'interactions'
$_{\Sigma_1}\mathcal{M}_{1\;\Sigma_2}$ and $_{\Sigma_2}\mathcal{M}_{2\;\Sigma_3}$ is
an 'interaction' $_{\Sigma_1}(\mathcal{M}_{1\;\Sigma_2}\mathcal{M}_2)_{\Sigma_3}$,
where $\mathcal{M}_{1\;\Sigma_2}\mathcal{M}_2$ is a result of gluing $\mathcal{M}_1$
and $\mathcal{M}_2$ along $\Sigma_2$.
\end{definition}
The trivial gluing along the empty boundary component is also admissible. For
instance we can glue a 'collapsion' of $\Sigma$ and the corresponding 'expansion' in
the trivial way. In this way we obtain an 'interaction'
$_{\Sigma}(\mathcal{M}\mathcal{M}^{op})_{\Sigma}$. If we glue the 'expansion' of
$\Sigma$ and the 'collapsion' of $\Sigma$ along $\Sigma$, then we obtain a class of
manifolds with empty boundaries.
\begin{example}
Classes of two dimensional surfaces with holes provide examples of string
interactions.
\end{example}
We wish to built the temporal support semigroup as an arbitrary sequence
\begin{equation}
\begin{array}{c}
X_{1}\stackrel{f_{1}}{\longrightarrow }X_{2}\stackrel{f_{2}}{\longrightarrow }\cdots
\stackrel{f_{n-1}}{\longrightarrow }X_{n} \label{seq}
\end{array}
\end{equation}
of objects and arrows of a directed graph $\mathfrak{C}$ indexed by a discrete time.
We wish to represent an 'interaction' $_{\Sigma_1}\mathcal{M}_{\Sigma_2}$ as an arrow
$X_1\stackrel{f}{\rightarrow}X_2$ of $\mathfrak{C}$. Obviously composable arrows
$X_{1}\stackrel{f_1}{\longrightarrow }X_{2}\stackrel{f_{2}}{\longrightarrow }X_3$
should represent the gluing
$_{\Sigma_1}(\mathcal{M}_{1\;\Sigma_2}\mathcal{M}_2)_{\Sigma_3}$. Two 'interactions'
$_{\Sigma_1}\mathcal{M}_{\Sigma_2}$ and $_{\Sigma'_1}\mathcal{M'}_{\Sigma'_2}$ should
be represented by the same arrow $X_1\stackrel{f}{\rightarrow}X_2$ if and only if both
'interactions' are 'parallel (simultaneous) in the time'.

Let us assume that the directed graph $\mathfrak{C}$ is an n-regular monoidal
category with duals. Let $X_{1}\stackrel{f_{1}}{\longrightarrow
}X_{2}\stackrel{f_{2}}{\longrightarrow }\cdots \stackrel{f_{n-1}}{\longrightarrow
}X_{n}$ be an n-regular cocycle. If there is an equivalence $\cong$ in
$\mathfrak{C}ob$ such that objects of the n-regular cocycle represent equivalence
classes of $\cong$ and arrows represent time consequnces, then we say that we have an
n-regular TQFT.

What means here the n-regularity? It is natural to assume that the opposite
$_{\Sigma_2}\mathcal{M}^{op}_{\Sigma_1}$ of $_{\Sigma_1}\mathcal{M}_{\Sigma_2}$ should
be representing by a reversed arrow $X_1\stackrel{f}{\leftarrow}X_2$.  The trivial
2-regularity is clear, it means that the time is reversible. We postulate that the
time is directed and always run further, never back, never stop. In other words, 'our
time' is not reversible in general, but it can be n-regular, where the regularity is
nontrivial.
\begin{example}
Let us consider for instance the 2-regular 'interactions'. Let
$$_{\Sigma_1}\mathcal{M}_{1\;\Sigma_2}\;\mbox{and}\;
 _{\Sigma_2}\mathcal{M}_{2\;\Sigma_1}$$ be two interactions, then
$_{\Sigma_1}(\mathcal{M}_{1\;\Sigma_2}\mathcal{M}_2)_{\Sigma_1}$ and
$_{\Sigma_2}(\mathcal{M}_{2\;\Sigma_1}\mathcal{M}_1)_{\Sigma_2}$ can be represented
as arrows $X_1\stackrel{f_1}{\rightarrow}X_2\stackrel{f_2}{\rightarrow}X_1$, and
$X_2\stackrel{f_2}{\rightarrow}X_1\stackrel{f_1}{\rightarrow}X_2$, respectively.
Interactions
$_{\Sigma_1}\mathcal{M}_{1\;\Sigma_2}\mathcal{M}_{2\;\Sigma_1}\mathcal{M}_{1\;\Sigma_2}$
and
$_{\Sigma_2}\mathcal{M}_{2\;\Sigma_1}\mathcal{M}_{1\;\Sigma_2}\mathcal{M}_{2\;\Sigma_1}$
should be represented by
$X_1\stackrel{f_1}{\rightarrow}X_2\stackrel{f_2}{\rightarrow}X_1
\stackrel{f_1}{\rightarrow}X_2$, and
$X_2\stackrel{f_2}{\rightarrow}X_1\stackrel{f_1}{\rightarrow}X_2
\stackrel{f_2}{\rightarrow}X_1$, respectively. Now the 2-regularity conditions are
clear.
\end{example}
Observe that the regularity concept can be useful for the construction of quantum
theory of the whole universe with nonreversible time evolution. In fact the
nontrivial n-regularity conditions  mean that all processes always go further, never
back, never stop, but are cyclically repeating after n-steps up to an equivalence.

\bigskip

\textbf{Acknowledgments}. One of the authors (S.D.) would like to thank Jerzy
Lukierski for kind hospitality at the University of Wroclaw, where this work was
begun, to Andrzej Borowiec,  Andrzej Frydryszak,
 and Cezary
Juszczak for valuable help during his stay in Wroclaw, also he is thankful to
Friedemann Brandt, Dimitry Leites, Volodymyr Lyubashenko for fruitful discussions at
the NATO ARW in Kiev and Fang Li for useful correspondence and rare reprints.

The paper is partially supported by the Polish KBN Grant No 5P03B05620.


\begin{thebibliography}{9}
\bibitem{gell} M. Gell-Mann and J. Hartle, Quantum mechanics in the light of quantum
cosmology, in Proceedings of the Third International Symphosium on the Foundations of
Quantum Mechanics in the light of New Technology, ed. by S. Kobayashi et al,
(Physical Society of Japan, Tokyo 1990), pp 321-343.
\bibitem{isham} C. J. Isham, Quantum Logic and the histories approach to quantum
theory, J. Math. Phys. \textbf{35}, 2157, (1994).
\bibitem{at} M. Atiyah, \textit{The Geometry and Physics of Knots},
Cambridge 1990.
\bibitem{bad} J. C. Baez and J. Dolan,  Higher--dimensional algebra
and topological quantum field theory, \textit{Journal of Mathematical Physics}
\textbf{36}, 6073. (1996).
\bibitem{saw} S. Sawin, Links, quantum groups and TQFT's,
{\it Bull. (New Series) AMS} {\bf 33}, 413 (1996).
\bibitem{ker} T. Kerler, Bridged links and tangle
presentation of cobordism categories, math.GT/980614 (1998).
\bibitem{cra} L. Crane, J. Math. Phys. {\bf 36}, 6180 (1995).
\bibitem{cran} L. Crane, Toplogical Field Theory as the Key to
Quantum Gravity, lecture presented to the conference on knot theory and quantm
gravity Riverside, California.

\bibitem{neu}  J.~von Neumann, \newblock "On regular rings," \newblock Proc.
Nat. Acad. Sci. USA \newblock {\bf  22}, 707 (1936).

\bibitem{pen1}  R.~Penrose, \newblock "A generalized inverse for matrices," %
\newblock Math. Proc. Cambridge Phil. Soc. \newblock {\bf  51}, 406 (1955).

\bibitem{rabson}  G.~Rabson,
\newblock {\it The Generalized Inverses in Set
Theory and Matrix Theory} \newblock (Amer. Math. Soc., Providence, 1969).

\bibitem{rao/mit}  C.~R. Rao and S.~K. Mitra,
\newblock {\it Generalized
Inverse of Matrices and its Application} \newblock (Wiley, New York, 1971).

\bibitem{nashed}  M.~Z. Nashed,
\newblock {\it Generalized Inverses and
Applications} \newblock (Academic Press, New York, 1976).

\bibitem{mun/pen}  W.~D. Munn and R.~Penrose, \newblock "Pseudoinverses in
semigroups," \newblock Math. Proc. Cambridge Phil. Soc. \newblock {\bf  57}, 247
(1961).

\bibitem{cli5}  A.~H. Cliford, \newblock "The fundamental representation of
a regular semigroup," \newblock Semigroup Forum \newblock {\bf  10}, 84 (1975/76).

\bibitem{lal1}  G.~Lallement, \newblock "Structure theorems for regular
semigroups," \newblock Semigroup Forum \newblock {\bf  4}, 95 (1972).

\bibitem{howie}  J.~M. Howie,
\newblock {\it An Introduction to Semigroup
Theory} \newblock (Academic Press, London, 1976).

\bibitem{lawson}  M.~V. Lawson,
\newblock {\it Inverse Semigroups: {T}he
Theory of Partial Symmetries} \newblock (World Sci., Singapore, 1998).

\bibitem{dup6}  S.~Duplij, \newblock "On semigroup nature of superconformal
symmetry," \newblock J.~Math. Phys. \newblock {\bf  32}, 2959 (1991).

\bibitem{dup11}  S.~Duplij, \newblock "Some abstract properties of
semigroups appearing in superconformal theories," \newblock Semigroup Forum %
\newblock {\bf  54}, 253 (1997).

\bibitem{duplij}  S. Duplij, \newblock {\it Semisupermanifolds and
semigroups}, \newblock (Krok, Kharkov, 2000).

\bibitem{dup18}  S.~Duplij, \newblock "On semi-supermanifolds," \newblock %
Pure Math. Appl. \newblock {\bf  9}, 1 (1998).

\bibitem{dup-hab}  S.~Duplij,
\newblock {\it Semigroup methods in supersymmetric theories of elementary
  particles} \newblock (Habilitation Thesis, Kharkov State University,
\texttt{\ math-ph/9910045}, Kharkov, 1999).

\bibitem{fangli3}  F.~Li, \newblock "Weak {H}opf algebras and new solutions
of {Y}ang-{B}axter equation," \newblock J. Algebra \newblock {\bf  208}, 72 (1998).

\bibitem{fangli2}  F.~Li, \newblock "Solutions of {Y}ang-{B}axter equation
in an endomorphism semigroup and quasi-(co)braided almost bialgebras,", Zhejiang
Univ. \textit{preprint}, Hangzhou, 1999.

\bibitem{nill} F, Nill, {\it Axioms for  week bialgebras}, {\tt q-alg/9805104}, 1998.

\bibitem{fangdup} S.~Duplij and F.~Li, \newblock "Regular solutions of quantum
Yang-Baxter equation from weak Hopf algebras", Czechoslovak Journal of Physics
51(12): 1306-1311; Dec (2001), math.QA/0105064, (2001).

\bibitem{fangli1} F. Li and S. Duplij, \newblock "Weak Hopf Algebras and Singular
Solutions of Quantum Yang-Baxter Equation", Commun. Math. Phys. 2002, V.225, N1,
191-217.

\bibitem{dav/rob}  D.~L. Davis and D.~W. Robinson, \newblock "Generalized
inverses of morphisms," \newblock Linear Algebra Appl. \newblock {\bf  5}, 329 (1972).

\bibitem{rob/cap}  J.~A.~G. Roberts and H.~W. Capel, \newblock "Area
preserving mappings that are not reversible," \newblock Phys. Lett. %
\newblock {\bf  A162}, 243 (1992).

\bibitem{ara}  A.~Arai, \newblock "Noninvertible {B}ogolyubov
transformations and instability of embedded eigenvalues," \newblock J. Math. Phys.
\newblock {\bf  32}, 1838 (1991).

\bibitem{mitchell}  B.~Mitchell, \newblock {\it Theory of Categories} %
\newblock (Academic Press, New York, 1965).

\bibitem{maclane1}  S.~MacLane,
\newblock {\it Categories for the Working
Mathematician}, The Second Edition, \newblock (Springer-Verlag, Berlin, 2000).

\bibitem{joy/str}  A.~Joyal and R.~Street, \newblock "Braided monoidal
categories,", Macquarie University \textit{preprint}, Mathematics Reports 86008,
North Ryde, New South Wales, 1986.

\bibitem{maj3}  S.~Majid, \newblock "Quasitriangular {H}opf algebras and {Y}%
ang-{B}axter equations," \newblock Int. J. Mod. Phys. \newblock {\bf  A5}, 1 (1990).

\end{thebibliography}
\end{document}